\documentclass{amsart}
\usepackage{subfigure}
\usepackage{graphicx}
\usepackage{hyperref}
\pdfoutput=1
\allowdisplaybreaks[4]
\numberwithin{equation}{section}

\newcommand{\ds}{\displaystyle}
\newcommand{\ts}{\textstyle}

\def\Mod{\operatorname{mod}}
\def\textcap{\,{\ts \bigcap}\,}

\def\RR{{\mathbb R}}
\def\CM{{\mathcal{M}}}
\def\bM{\bar{M}}

\def\da{\delta a}
\def\dth{\delta \theta}
\def\dze{\delta \zeta}

\newtheorem*{thmFen}{Theorem (Fenichel \cite{Fe71})}
\newtheorem{lem}{Lemma}[section]

\begin{document}

\title[Breakdown of Normal Hyperbolicity]{Breakdown of Normal Hyperbolicity for a Family of Invariant Manifolds with Generalized Lyapunov-Type Numbers Uniformly Bounded Below Their Critical Values}

\author{Dennis Guang Yang}
\address{Department of Mathematics, Cornell University, Ithaca, New York 14853}
\email{gy26@cornell.edu}


\date{December, 2009}



\begin{abstract}
We present three examples to illustrate that in the continuation of a family of normally hyperbolic $C^1$ manifolds, the normal hyperbolicity may break down as the continuation parameter approaches a critical value even though the corresponding generalized Lyapunov-type numbers remain uniformly bounded below their critical values throughout the process. In the first example, a $C^1$ manifold still exists at the critical parameter value, but it is no longer normally hyperbolic. In the other two examples, at the critical parameter value the family of $C^1$ manifolds converges to a nonsmooth invariant set, for which generalized Lyapunov-type numbers are undefined. 
\end{abstract}

\maketitle



\section{Introduction} \label{intro}

In normally hyperbolic invariant manifold theory, generalized Lyapunov-type numbers characterize the linearized dynamics along an invariant manifold and provide an effective way to determine the normal hyperbolicity of the invariant manifold. They are introduced by Fenichel in his seminal work on persistence and smoothness of normally hyperbolic invariant manifolds for flows \cite{Fe71}. Analogous concepts are also developed for maps by others (see, e.g., \cite{HiPuSh77}). Roughly speaking, an invariant manifold is normally hyperbolic and persists under small perturbations if the corresponding generalized Lyapunov-type numbers are less than some critical values. 

A particularly interesting problem is the continuation of a family of normally hyperbolic invariant manifolds with respect to some parameter until the breakdown of normal hyperbolicity at some critical parameter value. Intuitively, one might expect that near the breakdown of normal hyperbolicity, some of the generalized Lyapunov-type numbers of these manifolds would have to become arbitrarily close to their critical values. However, this is in fact not the case. In \cite{ChLi00} Chicone and Liu describe a scenario in which a family of normally hyperbolic limit cycles converges to a nonsmooth homoclinic loop and thus loses normal hyperbolicity even through the generalized Lyapunov-type numbers for the whole family can be made uniformly bounded away from their critical values. In addition, in \cite{HaLl07} Haro and de la Llave study numerical continuation of invariant tori for quasi-periodic perturbations of the standard map and the H\'{e}non map. They observe a situation where a family of normally hyperbolic invariant 1-tori cannot be continued further as their stable and unstable bundles converge locally when the continuation parameter approaches a critical value. However, the corresponding generalized Lyapunov-type numbers remain uniformly bounded below their critical values throughout the process.

Compared to the vast literature on the application of normally hyperbolic invariant manifolds, little attention has been paid to this counterintuitive behavior of generalized Lyapunov-type numbers. In fact, to the best knowledge of the author, the two aforementioned examples are by far the only references in the literature. The purpose of this paper is to provide further examples to illustrate this behavior of generalized Lyapunov-type numbers and to show that how much the generalized Lyapunov-type numbers are less than their critical values gives no information about how robust the normal hyperbolicity of an invariant manifold is. In the following sections, we first give an overview of generalized Lyapunov-type numbers and then discuss three examples involving invariant manifolds for systems of ODEs.

\section{An Overview of Generalized Lyapunov-Type Numbers} \label{overview}

Following the presentation of Fenichel \cite{Fe71}, we give the following definitions of generalized Lyapunov-type numbers. Consider a vector field defined on $\RR^n$:
\begin{equation} \label{fx}
\dot{x} = f(x)\,, 
\end{equation}
where $x \in \RR^n$, $F$ is $C^r$ for some $r \ge 1$, and $\dot{} = \tfrac{d}{dt}$. Let $\phi^t$ be the flow of (\ref{fx}). Suppose $\bM = M \cup \partial M \subset \RR^n$ is a $C^r$, compact, and connected manifold with boundary and overflowing invariant under the flow $\phi^t$. Let $TM$ denote the tangent bundle of $M$. In addition, let $N$ denote the normal bundle of $M$ with respect to the standard inner product on $\RR^n$. Then, we have the splitting
\begin{equation*}
T \RR^n |_{M} = TM \oplus N
\end{equation*}
with the associated orthogonal projection
\begin{equation*}
\Pi : T \RR^n |_{M} \rightarrow N \,.
\end{equation*}
For each $p \in M$, define the linear operators $A_t(p)$ and $B_t(p)$ constructed from the linearized flow $D \phi^{t}$ as follows:
\begin{equation} \label{def:AB}
\begin{split}
A_t(p) :={}& D \phi^{-t}(p) |_{T_p M} \,:\, T_p M \rightarrow T_{\phi^{-t}(p)} M \,, \\
B_t(p) :={}& \Pi_p  D \phi^t (\phi^{-t}(p)) |_{N_{\phi^{-t}(p)}} \,:\, N_{\phi^{-t}(p)} \rightarrow N_p \,.
\end{split}
\end{equation}
Then, the generalized Lyapunov-type numbers $\nu(p)$ and $\sigma(p)$ are defined as follows:
\begin{equation} \label{def:LN}
\begin{split}
\nu(p) :={}& \limsup_{t\rightarrow\infty} \| B_t(p) \|^{\frac{1 }{ t}} \,, \\
\sigma(p) :={}& \limsup_{t\rightarrow\infty} \frac{ \log{\| A_t(p)\|} }{ -\log{\| B_t(p)\|} } \,.
\end{split}
\end{equation}
The number $\nu(p)$ measures the exponential of the normal contraction rate under $D \phi^{t}(p)$, and the number $\sigma(p)$ compares the tangential and normal contraction rates under $D \phi^{t}(p)$. In addition, both numbers are constant along the trajectory $\phi^{-t}(p)$ for all $t \ge 0$, and $\nu(p) = \nu(p')$ and $\sigma(p) = \sigma(p')$ for any $p'$ in the $\alpha$-limit set of $p$. The proofs of these properties can be found in \cite{Fe71,Wi94}. In \cite{Fe71} Fenichel has also introduced other generalized Lyapunov-type numbers in the same spirit as the above construction but for a more general situation where the normal bundle of $M$ further splits into stable and unstable subbundles. However, the examples considered in this paper involve only the two defined by (\ref{def:LN}). 

If $\nu(p) < 1$ and $\sigma(p) < \tfrac{1 }{ r}$ for all $p \in M$, then the overflowing invariant manifold $\bM$ is said to be normally hyperbolic (or more precisely, $r$-normally hyperbolic) and $\bM$ persists under small perturbations according to the following theorem.

\begin{thmFen}
Let\/ $\dot{x} = f(x)$ be a\/ $C^r\!$ vector field on\/ $\RR^n$, $r \ge 1$. Let\/ $\bM = M \cup \partial M$ be a\/ $C^r\!$, compact, and connected manifold with boundary and overflowing invariant under the flow of\/ $\dot{x}=f(x)$. Suppose\/ $\nu(p) < 1$ and\/ $\sigma(p) < \tfrac{1 }{ r}$ for all\/ $p \in M$. Then for any\/ $C^r\!$ vector field\/ $f_{\epsilon}$ in a sufficiently small\/ $C^1\!$ neighborhood of\/ $f$, there is a manifold\/ $\bM_{\epsilon}$ overflowing invariant under the flow of\/ $\dot{x}=f_{\epsilon}(x)$ and\/ $C^r\!$ diffeomorphic to\/ $\bM$.
\end{thmFen}

When $M$ is a $C^r$, compact, and connected boundaryless manifold and invariant under the flow $\phi^t$, the above construction of generalized Lyapunov-type numbers using (\ref{def:AB}) and (\ref{def:LN}) is still applicable. Furthermore, suppose $\nu(p) < 1$ and $\sigma(p) < \tfrac{1 }{ r}$ for all $p \in M$. Then we can retain the notion of normal hyperbolicity and the accompanying persistence property even though $M$ is boundaryless. Specifically, we can append an auxiliary variable $\zeta \in \RR$ to (\ref{fx}) to form an enlarged system as follows:
\begin{equation} \label{fx-zeta}
\begin{split}
\dot{x} ={}& f(x)\,, \\
\dot\zeta ={}& \zeta^3 \,.
\end{split} 
\end{equation}
Then $\CM := M \times [-1, 1] \subset \RR^{n+1}$ is a $C^r$, compact, and connected manifold with boundary and overflowing invariant under the flow of (\ref{fx-zeta}). It is straightforward to verify that for any $(p,\zeta) \in \CM$, $\tilde\nu(p,\zeta) = \nu(p)$ and $\tilde\sigma(p,\zeta) = \sigma(p)$, where $\tilde\nu(p,\zeta)$ and $\tilde\sigma(p,\zeta)$ are the corresponding generalized Lyapunov-type numbers computed using the flow of (\ref{fx-zeta}). Then $\CM$ is $r$-normally hyperbolic with respect to the flow of (\ref{fx-zeta}), and it perturbs to a $C^r$ manifold $\CM_{\epsilon}$ that is overflowing invariant under the flow of
\begin{align*}
\dot{x} ={}& f_{\epsilon}(x)\,, \\
\dot\zeta ={}& \zeta^3 \,,
\end{align*} 
for any $C^r$ vector field $f_{\epsilon}$ that is sufficiently $C^1$-close to $f$. Clearly, the cross-section of $\CM_{\epsilon}$ at $\zeta = 0$ is $C^r$ diffeomorphic to $M$ and invariant under the flow of $\dot{x} = f_{\epsilon}(x)$. In the balance of this paper, we consider only boundaryless invariant manifolds. Furthermore, we take $r = 1$ and refer to a $C^1$, compact, and connected boundaryless manifold $M$ as being normally hyperbolic if $\nu(p) < 1$ and $\sigma(p) < 1$ for all $p \in M$.

\section{Example 1} \label{ex1}

Consider the following $1$-parameter family of $2$-dimensional systems:
\begin{equation} \label{sys1}
\begin{split}
\dot a ={}& -\left( \tfrac{1}{2} + \sin \theta \right) a \,, \\
\dot \theta ={}& \beta \,,
\end{split}
\end{equation}
where $a \in \RR$, $\theta \in \RR (\Mod 2\pi)$, and the parameter $\beta \in \RR$. Clearly, the circle 
\begin{equation*}
\Gamma := \big\{ (a, \theta) : a = 0,\, \theta \in \RR (\Mod 2\pi) \big\}
\end{equation*} 
is an invariant manifold for (\ref{sys1}) with any $\beta \in \RR$. In particular, $\Gamma$ is a periodic orbit if $\beta \neq 0$, and it is a set of fixed points if $\beta = 0$. Let us compute the generalized Lyapunov-type numbers for $\Gamma$.

Let $\phi^t_{\beta}$ be the flow of (\ref{sys1}) with the corresponding $\beta$. Note that for any $\beta \in \RR$, $\phi^t_{\beta}(p) = (0, \theta + \beta t)$ for any $p = (0, \theta) \in \Gamma$. Suppose $\beta \neq 0$. Then straightforward application of (\ref{def:AB}) and (\ref{def:LN}) yields that for any point $p$ in the periodic orbit $\Gamma$, $\sigma_{\beta}(p) = 0$ and 
\begin{equation*}
\nu_{\beta}(p) = e^{ - \frac{1}{2\pi} \int^{2\pi}_0 \left( \frac{1}{2} + \sin \theta \right) d\theta } = e^{-\frac{1}{2}} \,,
\end{equation*}
which is also the Floquet multiplier of any periodic solution along $\Gamma$. It follows that $\Gamma$ is normally hyperbolic with respect to the flow of (\ref{sys1}) with any $\beta \neq 0$. However, for (\ref{sys1}) with $\beta = 0$, $\Gamma$ is not normally hyperbolic since 
\begin{equation*}
\nu_0(p) = e^{-\left( \frac{1}{2} + \sin \theta \right)} \ge 1
\end{equation*}
for any $p = (0, \theta) \in \Gamma$ with $\theta \in [\frac{7}{6}\pi, \frac{11}{6}\pi]$. Therefore, $\Gamma$ loses normal hyperbolicity at $\beta = 0$ even though $\Gamma$ is normally hyperbolic with the generalized Lyapunov-type numbers $\nu_{\beta}(p) = e^{-\frac{1}{2}} < 1$ and $\sigma_{\beta}(p) = 0 < 1$ for any $p \in \Gamma$ and for any $\beta \neq 0$.

\section{Example 2} \label{ex2}

As mentioned earlier, in \cite{ChLi00} Chicone and Liu describe a scenario in which a family of normally hyperbolic limit cycles converges to a nonsmooth homoclinic loop while the generalized Lyapunov-type numbers for the whole family can be made uniformly bounded away from the critical value $1$. However, they have not provided an explicit formulation of the system that admits such a dynamical feature. Below we present a simple example based on the mechanism introduced in \cite{ChLi00}.

\begin{figure}
  \centering
    \subfigure[$-0.1 \le c < c_1$]{\label{fig:ChLi:1}\includegraphics[scale=1]{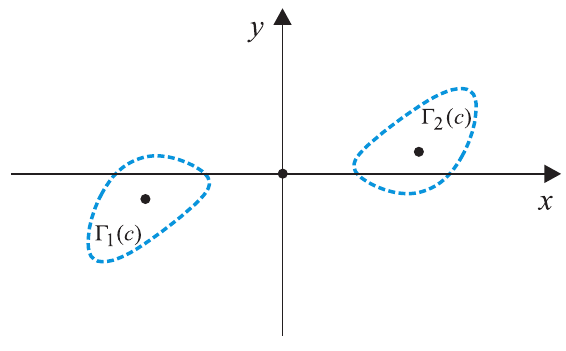}} \hspace{0.1in}
    \subfigure[$c=c_1 \simeq -0.060959$]{\label{fig:ChLi:2}\includegraphics[scale=1]{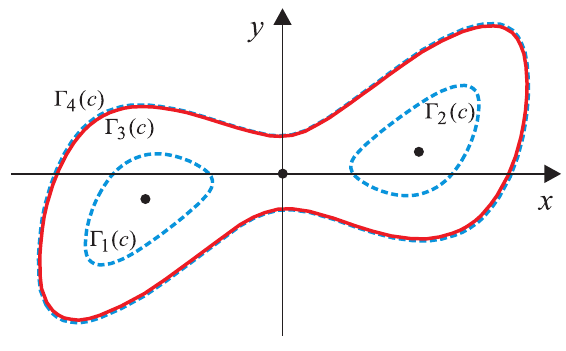}} \\
    [0.7em]
    \subfigure[$c_1 < c < c_2$]{\label{fig:ChLi:3}\includegraphics[scale=1]{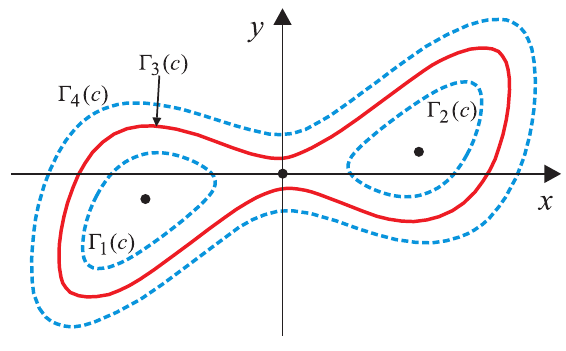}} \hspace{0.1in}
    \subfigure[$c=c_2 \simeq -0.060932$]{\label{fig:ChLi:4}\includegraphics[scale=1]{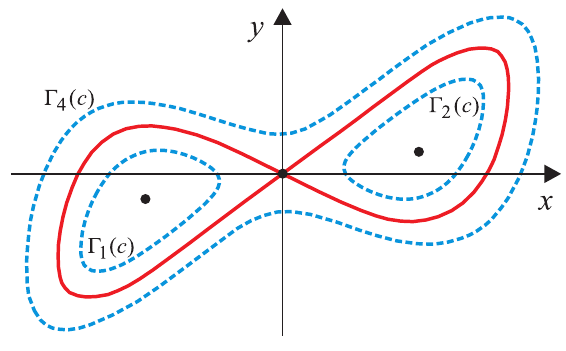}} \\
    [0.7em]
    \subfigure[$c_2 < c < c_3$]{\label{fig:ChLi:5}\includegraphics[scale=1]{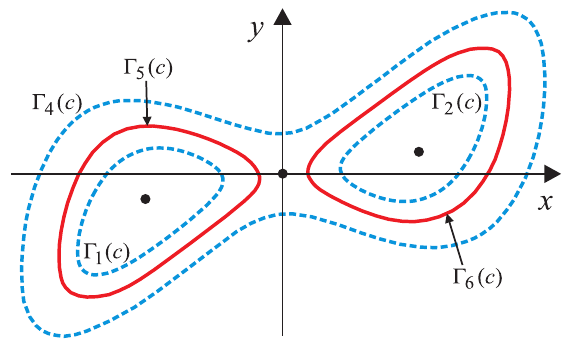}} \hspace{0.1in}
    \subfigure[$c=c_3 \simeq -0.060903$]{\label{fig:ChLi:6}\includegraphics[scale=1]{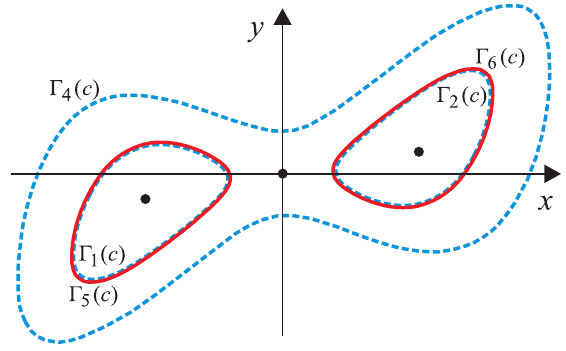}} \\
    [0.7em]
    \subfigure[$c_3 < c \le 0$]{\label{fig:ChLi:7}\includegraphics[scale=1]{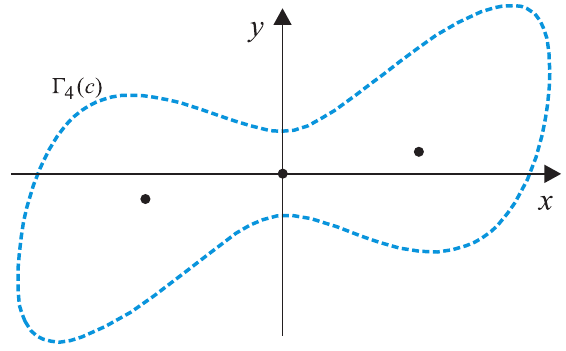}}
  \caption{Sketches of the phase portrait of (\ref{sys:ChLi}) for various values of $c$. Black dots denote fixed points, blue dotted curves are unstable limit cycles, and red solid curves are stable limit cycles or stable homoclinic orbits (in (d)). Note that the variations of the limit cycles are largely exaggerated. All the limit cycles are actually contained in a small neighborhood of the homoclinic loop.}
  \label{fig:ChLi}
\end{figure}

Consider a $1$-parameter family of planar systems:
\begin{equation} \label{sys:ChLi}
\begin{split}
\dot{x} ={}& -\tfrac{1}{2} x + y \,, \\
\dot{y} ={}& x - x^3 + y \left( \tfrac{1}{2} y^2 - \left( \tfrac{1}{2} x^2 - \tfrac{1}{4} x^4 \right) - c \right) \,,
\end{split}
\end{equation}
where the parameter $c \in [-0.1, 0]$. It is clear that the origin $(0,0)$ is a saddle point for any $c \in [-0.1, 0]$. There are three critical parameter values: $c_1 \simeq -0.060959$, $c_2 \simeq -0.060932$, and $c_3 \simeq -0.060903$, at which bifurcations occur. For $-0.1 \le c < c_1$, (\ref{sys:ChLi}) has two unstable limit cycles $\Gamma_1(c)$ and $\Gamma_2(c)$ as sketched in Figure \ref{fig:ChLi:1}. At $c=c_1$, a saddle-node bifurcation of limit cycles creates a new pair of limit cycles $\Gamma_3(c)$ (stable) and $\Gamma_4(c)$ (unstable) as shown in Figure \ref{fig:ChLi:2}. For $c_1 < c < c_2$, $\Gamma_3(c)$ shrinks in size as $c$ increases (Figure \ref{fig:ChLi:3}) and eventually converges to a figure-$8$ homoclinic loop of the saddle point $(0,0)$ in the limit $c \rightarrow c_2^{-}$. Note that the homoclinic loop exists only when $c=c_2$, at which a homoclinic bifurcation occurs (Figure \ref{fig:ChLi:4}). For $c_2 < c < c_3$, $\Gamma_3(c)$ bifurcates into two smaller stable limit cycles $\Gamma_5(c)$ and $\Gamma_6(c)$ as shown in Figure \ref{fig:ChLi:5}. At $c = c_3$, $\Gamma_5(c)$ and $\Gamma_6(c)$ collide with $\Gamma_1(c)$ and $\Gamma_2(c)$, respectively, in two simultaneous saddle-node bifurcations of limit cycles (Figure \ref{fig:ChLi:6}). Finally, only the unstable limit cycle $\Gamma_4(c)$ remains for $c_3 < c \le 0$, as displayed in Figure \ref{fig:ChLi:7}.

Now consider the continuation of the stable limit cycle $\Gamma_3(c)$ for $c \in [c_1 + \delta, c_2)$ with $\delta$ being some small positive constant such that $c_1 + \delta < c_2$. Since $\Gamma_3(c)$ is a limit cycle, the generalized Lyapunov-type number $\sigma_c(p) = 0$ for all $p \in \Gamma_3(c)$. Furthermore, the generalized Lyapunov-type number $\nu_c(p) = \lambda_c$ for all $p \in \Gamma_3(c)$ with $\lambda_c$ being the Floquet multiplier of any periodic solution $\gamma(\cdot,c) : \RR \rightarrow \Gamma_3(c)$ along $\Gamma_3(c)$. For the planar system (\ref{sys:ChLi}), $\lambda_c$ can be computed by the following formula (see, e.g., \cite{Ch06})
\begin{equation*}
\lambda_c = e^{ \frac{1}{T(c)} \int^{T(c)}_0 \operatorname{div}\! f(\gamma(t,c),c)\, dt } \,,
\end{equation*}
where $f((x,y),c)$ denotes the right-hand side of (\ref{sys:ChLi}) and $T(c)$ is the period of $\gamma(t,c)$. Note that
\begin{equation*}
\lim_{c\rightarrow c_2^{-}} \lambda_c = e^{\operatorname{div}\! f((0,0),c)} = e^{-\frac{1}{2}-c_2} < e^{-0.4} \,.
\end{equation*}
Thus, if $c_1 + \delta$ is chosen sufficiently close to $c_2$, the generalized Lyapunov-type numbers $\nu_c(p)$ and $\sigma_c(p)$ are uniformly bounded below $1$ for any $p \in \Gamma_3(c)$ with any $c \in [c_1 + \delta, c_2)$. However, the continuation of $\Gamma_3(c)$ has to cease at $c = c_2$.

\section{Example 3} \label{ex3}

In this example, we consider the continuation of a normally hyperbolic invariant torus in a $1$-parameter family of $3$-dimensional systems. During the continuation process, the torus continuously deforms. However, it contains the same $\alpha$-limit set. Consequently, the generalized Lyapunov-type numbers of the invariant torus remain constant for all the parameter values for which the torus exists. 

Consider the following $1$-parameter family of systems:
\begin{equation} \label{BF}
\begin{split}
\dot a ={}& -a + \beta^2 \sin^2\! \zeta \sin \theta \,, \\
\dot \zeta ={}& \sin^2\! \zeta \,, \\
\dot \theta ={}& a \,,
\end{split}
\end{equation}
where $a \in \RR$, $\zeta \in \RR (\Mod{\pi})$, $\theta \in \RR (\Mod{2\pi})$, and the parameter $\beta \in [0,1]$. Clearly, (\ref{BF}) with $\beta = 0$ has a unique invariant torus
\begin{equation*} 
T_0 = \big\{ (a, \zeta, \theta) : a = 0,\, \zeta \in \RR (\Mod{\pi}),\, \theta \in \RR (\Mod{2\pi}) \big\} \,.
\end{equation*}
In what follows, we will study the continuation of $T_0$ for $\beta > 0$. 

It is evident that for any $\beta \in [0,1]$, the set of fixed points of (\ref{BF}) forms an invariant circle
\begin{equation*}
\Gamma = \big\{ (a, \zeta, \theta) : a = 0,\, \zeta = 0,\, \theta \in \RR (\Mod{2\pi}) \big\} \,.
\end{equation*}
Let $\phi_{\beta}^t$ be the flow of (\ref{BF}) with the corresponding $\beta$. For any $(0, \zeta, \theta) \in T_0$, we have 
\begin{equation*}
\lim_{t \rightarrow -\infty} \phi_0^t (0, \zeta, \theta) = (0, 0, \theta) \in \Gamma \,.
\end{equation*}
Thus, to compute the generalized Lyapunov-type numbers for points in $T_0$, we only need to consider points in $\Gamma$ since $\nu_0(0,\zeta,\theta) = \nu_0(0,0,\theta)$ and $\sigma_0(0,\zeta,\theta) = \sigma_0(0,0,\theta)$ for any $(0,\zeta,\theta) \in T_0$. 

For any $(0, 0, \theta) \in \Gamma$, the linear variational equation along the solution trajectory $\phi_0^t (0, 0, \theta) \equiv (0, 0, \theta)$ is given by
\begin{equation} \label{BF:var}
\begin{pmatrix}
\dot{\da} \\ \dot{\dze} \\ \dot{\dth}
\end{pmatrix}
= 
\begin{pmatrix}
-1 & 0 & 0 \\
 0 & 0 & 0 \\
 1 & 0 & 0 
\end{pmatrix}
\begin{pmatrix}
\da \\ \dze \\ \dth
\end{pmatrix} .
\end{equation}
By integrating (\ref{BF:var}), we obtain
\begin{equation*} 
D\phi_0^t (0, 0, \theta)
= 
\begin{pmatrix}
e^{-t} & 0 & 0 \\
0 & 1 & 0 \\
1-e^{-t} & 0 & 1 
\end{pmatrix} .
\end{equation*}
Notice that at any $(0, 0, \theta) \in \Gamma$, the tangent space to $T_0$ has a basis $\big\{ \tfrac{\partial}{\partial \zeta}, \tfrac{\partial}{\partial \theta} \big\}$ and the normal space to $T_0$ has a basis $\big\{ \tfrac{\partial}{\partial a} \big\}$. Then straightforward computation using the definitions (\ref{def:AB}) and (\ref{def:LN}) yields $\nu_0(0,0,\theta) = e^{-1}$ and $\sigma_0(0,0,\theta) = 0$ for all $(0, 0, \theta) \in \Gamma$, which verifies the normal hyperbolicity of $T_0$ with respect to $\phi_0^t$. Therefore, for some $\beta^* > 0$, there exists a continuous family of $C^1$ tori $\big\{ T_{\beta} : \beta \in [0, \beta^*) \big\}$  such that each $T_{\beta}$ is invariant under the flow $\phi_{\beta}^t$ generated by $(\ref{BF})$ with the corresponding $\beta$.

For each $\beta \in [0, 1]$, let $\Lambda_{\beta}$ be the largest invariant set under the flow $\phi_{\beta}^t$ inside the domain $U$ defined as follows:
\begin{equation*}
U := \big\{ (a, \zeta, \theta) : |a| \le 2,\, \zeta \in \RR (\Mod{\pi}),\, \theta \in \RR (\Mod{2\pi}) \big\} \,.
\end{equation*}
In the following lemma, we state two properties of $\Lambda_{\beta}$ that are useful in our subsequent analysis of the continuation of $T_{\beta}$ for $\beta > 0$. 

\begin{lem} \label{lem1}
For each\/ $\beta \in [0, 1]$, $\Lambda_{\beta}$ has the following properties:
\begin{enumerate}
\item If\/ $(a, \zeta, \theta) \in \Lambda_{\beta}$ and\/ $\zeta \in [0, \frac{\pi}{2}]$, then\/ $|a| \le \sin^2\! \zeta$.
\item For each\/ $(a, \zeta, \theta) \in \Lambda_{\beta}$ with\/ $\zeta \in [0, \pi)$, there exists a\/ $\theta^* \in \RR (\Mod{2\pi})$ such that\/ $\ds{\lim_{t \rightarrow -\infty}} \phi_{\beta}^t (a, \zeta, \theta) = (0, 0, \theta^*) \in \Gamma$.
\end{enumerate}
\end{lem}

\begin{proof}
To prove property (1), we consider the time-reversed trajectories $\phi_{\beta}^{-\tau} (a, \zeta, \theta)$ for $\tau \ge 0$. Clearly, if $\zeta \in [0, \tfrac{\pi}{2}]$, then the $\zeta$-component of $\phi_{\beta}^{-\tau} (a, \zeta, \theta)$ is contained inside the interval $[0, \tfrac{\pi}{2}]$ for all $\tau \ge 0$. Let $u := a - \sin^2\! \zeta$. Then at any $(a, \zeta, \theta)$ with $\zeta \in [0, \tfrac{\pi}{2}]$ and $a - \sin^2\! \zeta > 0$, we have that under the flow $\phi_{\beta}^{-\tau}$,
\begin{align*}
\frac{du}{d\tau} ={}& -( -a + \beta^2 \sin^2\! \zeta \sin \theta - 2\sin^3\!\zeta \cos\zeta ) \\
={}& u + \sin^2\! \zeta - \beta^2 \sin^2\! \zeta \sin \theta + 2\sin^3\!\zeta \cos\zeta \\
\ge{}& u > 0 \,.
\end{align*}
It follows that $u(\tau) \rightarrow \infty$ at least exponentially along $\phi_{\beta}^{-\tau} (u_0 + \sin^2\!\zeta, \zeta, \theta)$ as $\tau \rightarrow \infty$ if $u(0) = u_0 > 0$ and $\zeta \in [0, \frac{\pi}{2}]$. Next, let $v := a + \sin^2\! \zeta$. At any $(a, \zeta, \theta)$ with $\zeta \in [0, \tfrac{\pi}{2}]$ and $a + \sin^2\! \zeta < 0$, we have that under the flow $\phi_{\beta}^{-\tau}$,
\begin{align*}
\frac{dv}{d\tau} ={}& -( -a + \beta^2 \sin^2\! \zeta \sin \theta + 2\sin^3\!\zeta \cos\zeta ) \\
={}& v - \sin^2\! \zeta - \beta^2 \sin^2\! \zeta \sin \theta - 2\sin^3\!\zeta \cos\zeta \\
\le{}& v < 0 \,.
\end{align*}
Then $v(\tau) \rightarrow -\infty$ at least exponentially along $\phi_{\beta}^{-\tau} (v_0 - \sin^2\!\zeta, \zeta, \theta)$ as $\tau \rightarrow \infty$ if $v(0) = v_0 < 0$ and $\zeta \in [0, \frac{\pi}{2}]$. Therefore, for $(a, \zeta, \theta) \in \Lambda_{\beta}$ with $\zeta \in [0, \frac{\pi}{2}]$, the boundedness of the $a$-component of $\phi_{\beta}^{-\tau} (a, \zeta, \theta)$ for $\tau \rightarrow \infty$ makes it necessary that $|a| \le \sin^2\! \zeta$.

Recall that $\dot \zeta = \sin^2\! \zeta$. Consider any $(a, \zeta_0, \theta) \in \Lambda_{\beta}$ with $\zeta_0 \in (0, \pi)$, and reparametrize the negative semi-orbit originated at $(a, \zeta_0, \theta)$ by $\zeta \in (0, \zeta_0]$, i.e., 
\begin{equation*}
\big\{ (a_{\beta}(\zeta), \zeta, \theta_{\beta}(\zeta)) : \zeta \in (0, \zeta_0] \big\} := \big\{ \phi_{\beta}^{-\tau} (a, \zeta_0, \theta) : \tau \ge 0 \big\} \,.
\end{equation*}
Then property (1) implies that $|a_{\beta}(\zeta)| \le \sin^2\! \zeta$ for any $0 < \zeta \le \tfrac{\pi}{2}$. It follows that the limit
\begin{equation*}
\lim_{\zeta \rightarrow 0^+} \theta_{\beta}(\zeta) = \theta_{\beta}(\zeta_0) + \int_{\zeta_0}^{0} \frac{a_{\beta}(z)}{\sin^2\! z} dz
\end{equation*}
exists. Furthermore, if $(a, \zeta, \theta) \in \Lambda_{\beta}$ with $\zeta = 0$, then $a = 0$ and hence $\phi_{\beta}^t (a, 0, \theta) \equiv (0, 0, \theta)$ for all $t \in \RR$. Thus, property (2) is true.
\end{proof}

Now we take an arbitrary $\beta \in [0, \beta^*)$ and compute the generalized Lyapunov-type numbers for each point in $T_{\beta} \subseteq \Lambda_{\beta}$. In view of property (2) in Lemma \ref{lem1} and the obvious fact that $\Gamma \subset T_{\beta}$, this task can again be reduced to computing the generalized Lyapunov-type numbers for points in $\Gamma$ only. Note that for any $\beta > 0$ and any $(0, 0, \theta) \in \Gamma$, the linear variational equation of (\ref{BF}) along the solution trajectory $\phi_{\beta}^t (0, 0, \theta) \equiv (0, 0, \theta)$ is still the same as (\ref{BF:var}). Furthermore, property (1) implies that $T_{\beta}$ is tangent to $T_0$ along $\Gamma$. Thus the tangent space of $T_{\beta}$ at any point in $\Gamma$ coincides with that of $T_0$. Therefore, for any $\beta \in [0, \beta^*)$, we still have $\nu_{\beta}(0,0,\theta) = e^{-1}$ and $\sigma_{\beta}(0,0,\theta) = 0$ for all $(0, 0, \theta) \in \Gamma$. It follows that $\nu_{\beta}(a,\zeta,\theta) = e^{-1}$ and $\sigma_{\beta}(a,\zeta,\theta) = 0$ for any $(a,\zeta,\theta) \in T_{\beta}$ as long as $T_{\beta}$ exists.

\begin{figure}
  \centering
    \subfigure[Full view]{ \label{fig:Ls} \includegraphics[scale=1]{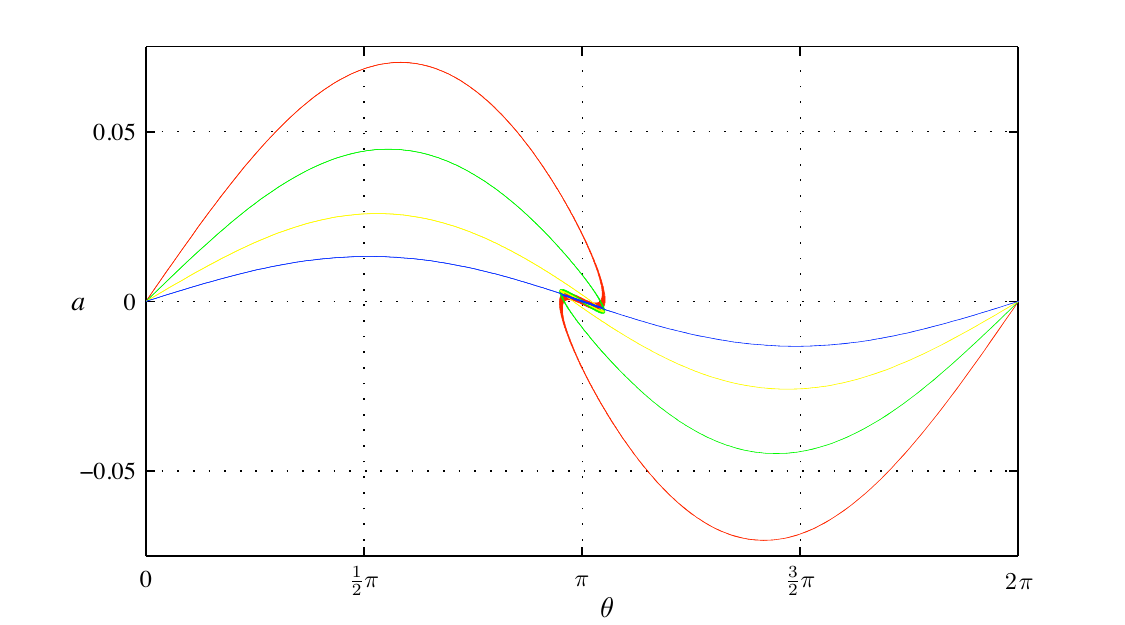} } \\
    [1em]
    \subfigure[Zoom view near $(a, \theta) = (0, \pi)$]{ \label{fig:Ls_zoom}\includegraphics[scale=1]{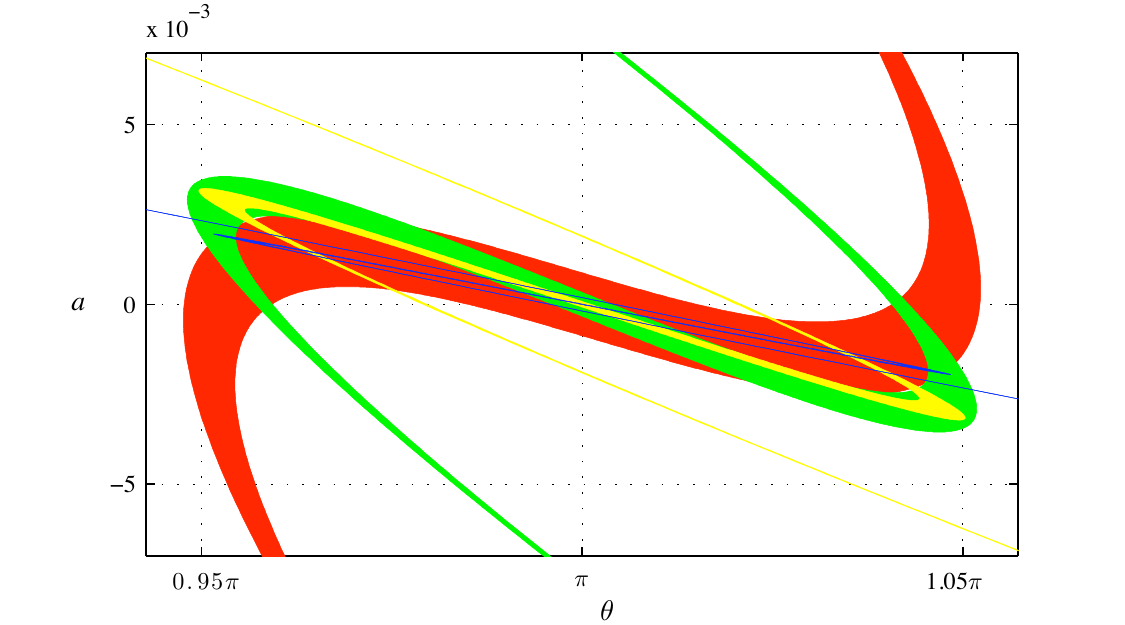} }
    \caption{The images of $\phi_1^t(L)$ at $\zeta = 0.89\pi$, $0.9\pi$, $0.91\pi$, and $0.92\pi$, plotted in red, green, yellow, and blue, respectively.}
  \label{fig:Lsecs}
\end{figure}

However, the system (\ref{BF}) with $\beta = 1$ has no $C^1$ invariant torus that can be continued from $T_{\beta}$. To see this, we consider the forward-time images of the set $L$ defined as follows:
\begin{equation*}
L := \big\{ (a, \zeta, \theta) : |a| \le \sin^2\! \tfrac{\pi}{20} + 10,\, \zeta = \tfrac{\pi}{20},\, \theta \in \RR (\Mod{2\pi}) \big\} \,.
\end{equation*}
Since $\dot \zeta = \sin^2\! \zeta$, by numerically integrating (\ref{BF}) with $\beta = 1$ starting at points on the boundary of $L$, we obtain a series of snapshots of $\phi_1^t(L)$ at various sections 
\begin{equation*}
\Sigma_{\zeta} = \big\{ (a, \zeta, \theta) : a \in \RR,\, \theta \in \RR (\Mod{2\pi}) \big\}
\end{equation*}
with monotonically increasing $\zeta$ converging to $\pi$ from below as $t \rightarrow \infty$. We find that a fold develops on $\phi_1^t(L)$ in a neighborhood of $(a, \theta) = (0, \pi)$ in the limit $t \rightarrow \infty$ as illustrated by the snapshots of $\phi_1^t(L)$ taken at $\Sigma_{\zeta}$ with $\zeta = 0.89\pi$, $0.9\pi$, $0.91\pi$, and $0.92\pi$ shown in Figure \ref{fig:Lsecs}. Suppose the family of $C^1$ tori $\big\{ T_{\beta} : \beta \in [0, \beta^*) \big\}$ can be continued up to $\beta = 1$. Then the $C^1$ torus $T_1$ must contain $\Gamma$, and $T_1 \textcap \Sigma_{\frac{\pi}{20}}$ must be contained in $L$ by property (1) in Lemma \ref{lem1}. It follows that $T_1 \textcap \Sigma_{\frac{\pi}{20}}$ is a circle winding around the ``cylinder'' $L$ exactly once. However, as $\zeta \rightarrow \pi^-$, the cross-section $T_1 \textcap \Sigma_{\zeta}$, which is contained in $\phi_1^t(L)$ taken at $\Sigma_{\zeta}$, has to fold as $\phi_1^t(L)$ does in the limit $t \rightarrow \infty$. This immediately contradicts the assumption that $T_1$ is a $C^1$ torus.

Let $\beta_c$ be the infimum of the set of $\beta$ with which (\ref{BF}) has no $C^1$ invariant torus that can be continued from $T_{\beta}$. The preceding analysis shows that $\beta_c \le 1$. In addition, it follows directly from the definition of $\beta_c$ that the family of $C^1$ tori $T_{\beta}$ exists for all $\beta \in [0, \beta_c)$. However, if (\ref{BF}) with $\beta = \beta_c$ has a $C^1$ invariant torus $T_{\beta_c}$, then $\nu_{\beta_c}(a,\zeta,\theta) = e^{-1}$ and $\sigma_{\beta_c}(a,\zeta,\theta) = 0$ for any $(a,\zeta,\theta) \in T_{\beta_c}$ since our previous computation of the generalized Lyapunov-type numbers for points in $T_{\beta}$ is still applicable. Consequently, $T_{\beta_c}$ is normally hyperbolic, and there exists a sufficiently small $\delta > 0$ such that the family of $C^1$ tori $T_{\beta}$ exists for all $\beta \in [0, \beta_c + \delta)$, which contradicts the definition of $\beta_c$. Thus, (\ref{BF}) with $\beta = \beta_c$ has no $C^1$ invariant torus that can be continued from $T_{\beta}$ although $T_{\beta}$ is normally hyperbolic and the corresponding generalized Lyapunov-type numbers are constant below $1$ for all $\beta \in [0, \beta_c)$.

\begin{figure}
  \centering
    \subfigure[$\beta = 0.65$]{\label{fig:RT1}\includegraphics[scale=1]{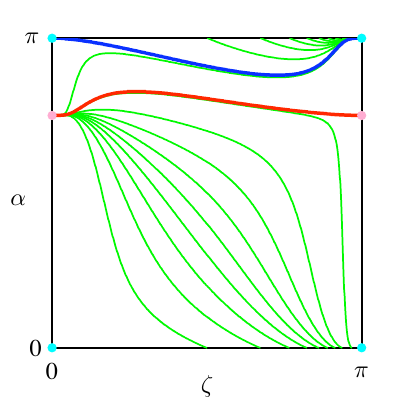}}
    \subfigure[$\beta = \beta_c \simeq 0.815$]{\label{fig:RT2}\includegraphics[scale=1]{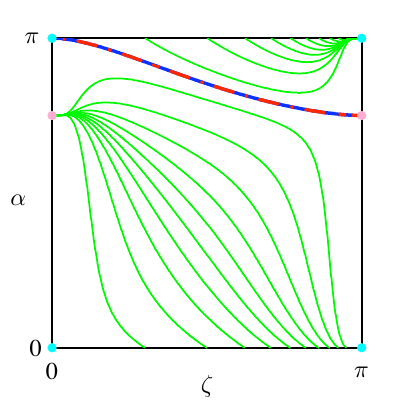}} 
    \subfigure[$\beta = 1$]{\label{fig:RT3}\includegraphics[scale=1]{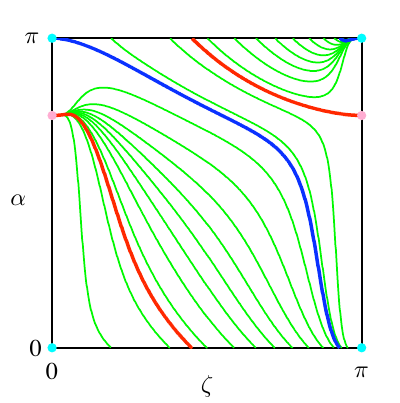}} 
  \caption{The phase portraits of (\ref{RT}) with $\beta = 0.65$, $\beta = \beta_c \simeq 0.815$, and $\beta = 1$. Note that $\zeta = 0$ and $\zeta = \pi$ are identified and $\alpha = 0$ and $\alpha = \pi$ are identified. The light-blue dots denote the fixed point $(\alpha, \zeta) = (0, 0)$, and the pink dots denote the fixed point $(\alpha, \zeta) = (\frac{3}{4}\pi, 0)$. In (a) and (c), the homoclinic orbits of $(\alpha, \zeta) = (0, 0)$ and $(\alpha, \zeta) = (\frac{3}{4}\pi, 0)$ are plotted in blue and red, respectively. In (b), the two homoclinic orbits ``collide'' and are replaced by a heteroclinic orbit shown in mixed blue and red.}
  \label{fig:RT}
\end{figure}

Further analysis reveals that the breakdown of $T_{\beta}$ in the limit $\beta \rightarrow \beta_c^-$ is due to the rotation of the normal bundle of the orbit
\begin{equation*}
E := \big\{ (a, \zeta, \theta) : a = 0,\, \zeta \in (0, \pi),\, \theta = \pi \big\} \,,
\end{equation*}
which is an invariant subset of $\Lambda_{\beta}$ for any $\beta \in [0, 1]$. Specifically, linearizing (\ref{BF}) about $E$ gives 
\begin{align*}
\dot \da ={}& -\da - \beta^2 \sin^2\! \zeta\, \dth \,, \\
\dot \zeta ={}& \sin^2\! \zeta \,, \\
\dot \dth ={}& \da \,,
\end{align*}
where $\da \in \RR$ and $\dth \in \RR$ are the variations in $a$ and $\theta$, respectively. Next, with the substitutions $\da = \rho \sin \alpha$ and $\dth = \rho \cos \alpha$, we obtain a system defined on the torus $\big\{ (\alpha, \zeta) : \alpha \in \RR (\Mod{\pi}),\, \zeta \in \RR (\Mod{\pi}) \big\}$ as follows:
\begin{equation} \label{RT}
\begin{split}
\dot \alpha ={}& - \cos \alpha \sin \alpha - \sin^2\! \alpha - \beta^2 \sin^2\! \zeta \cos^2\! \alpha \,, \\
\dot \zeta = {}& \sin^2\! \zeta \,,
\end{split}
\end{equation}
which describes the rotational dynamics in the normal bundle of $E$ under the linearized flow $D \phi_{\beta}^t (0, \zeta, \pi)$. For any $\beta \in [0,1]$, (\ref{RT}) has two fixed points $(\alpha, \zeta) = (0, 0)$ and $(\alpha, \zeta) = (\frac{3}{4}\pi, 0)$. For $\beta \in [0, \beta_c)$, each fixed point has a homoclinic orbit, along which the numbers of rotation in $\alpha$ and $\zeta$ are $0$ and $1$, respectively (see Figure \ref{fig:RT1} for the case with $\beta = 0.65$). As $\beta$ approaches $\beta_c$ from below, the two homoclinic orbits converge pointwise to each other. At $\beta = \beta_c \simeq 0.815$, the two homoclinic orbits ``collide'' and are replaced by a heteroclinic orbit (see Figure \ref{fig:RT2}). For $\beta \in (\beta_c, 1]$, the two fixed points are again connected to themselves by homoclinic orbits, but now the numbers of rotation in $\alpha$ and $\zeta$ are $-1$ and $1$, respectively, along both homoclinic orbits (see Figure \ref{fig:RT3} for the case with $\beta = 1$).

\begin{figure}
\centering
\includegraphics[scale=1]{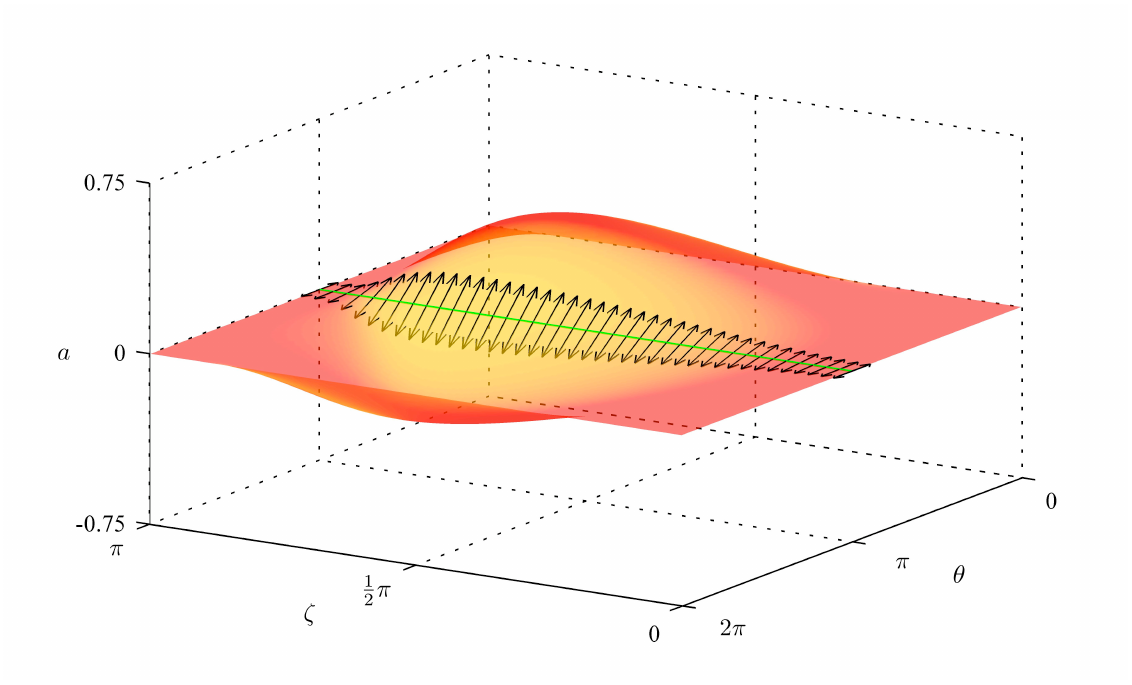}
\caption{$\Lambda_{\beta}$ with $\beta = 0.65$. The black arrows denote unit vectors in the fibers of $\mathcal{N}_{\beta}$ at various locations along the orbit $E$, which is shown in green. For any $\beta \in [0, \beta_c)$, $T_{\beta} = \Lambda_{\beta}$, and $\mathcal{N}_{\beta}$ is tangent to $T_{\beta}$ along the orbit $E$.}
\label{fig:smooth}
\end{figure}

\begin{figure}
\centering
\includegraphics[scale=1]{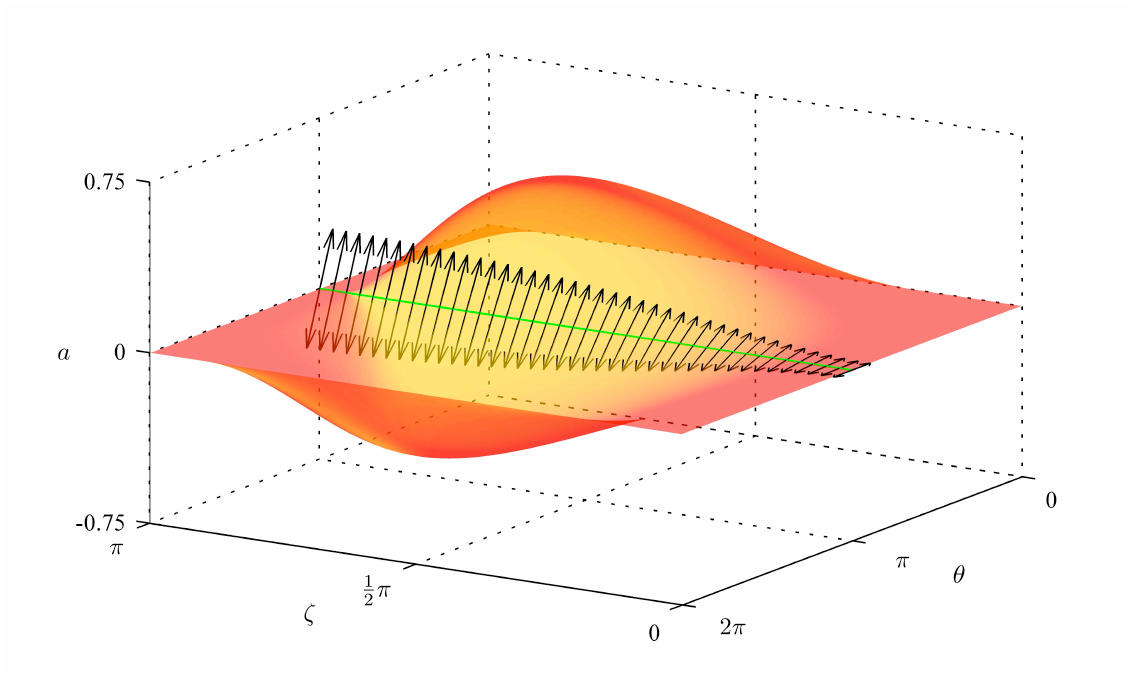}
\caption{$\Lambda_{\beta}$ with $\beta = \beta_c \simeq 0.815$. The black arrows denote unit vectors in the fibers of $\mathcal{N}_{\beta_c}$ at various locations along the orbit $E$, which is shown in green. $\Lambda_{\beta_c}$ is a nonsmooth topological torus.}
\label{fig:nonsmooth}
\end{figure}

\begin{figure}
  \centering
    \subfigure[Full view]{ \label{fig:fold} \includegraphics[scale=1]{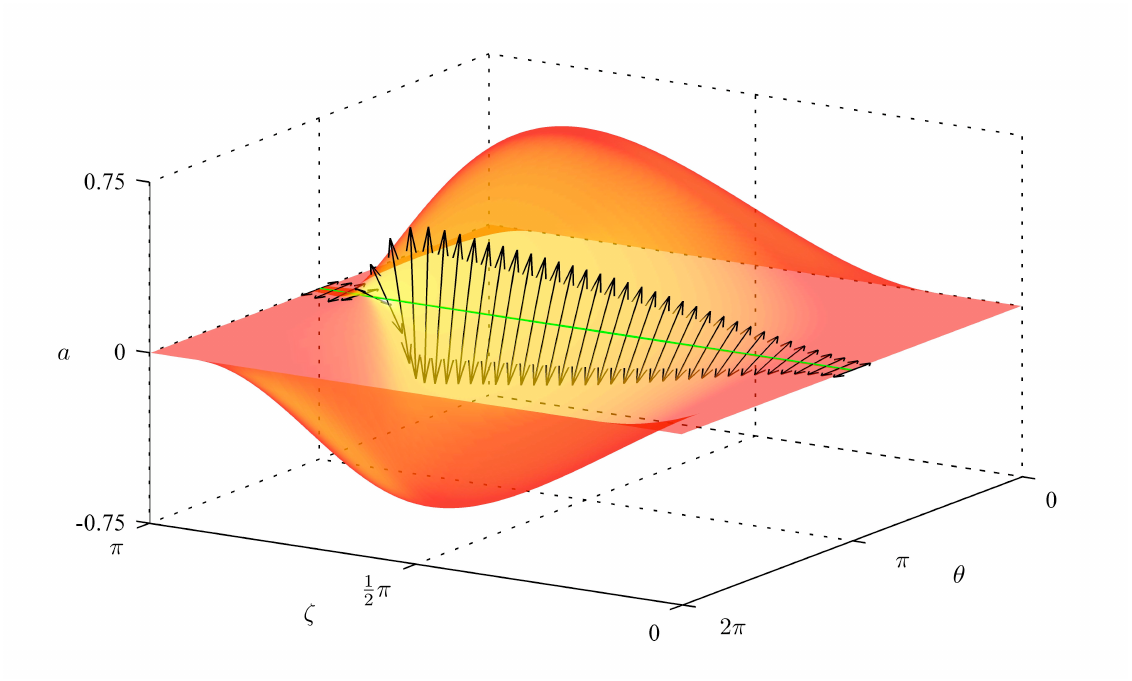} } \\
    \subfigure[Zoom view near $(a, \zeta, \theta) = (0, \pi, \pi)$]{ \label{fig:fold_zoom}\includegraphics[scale=1]{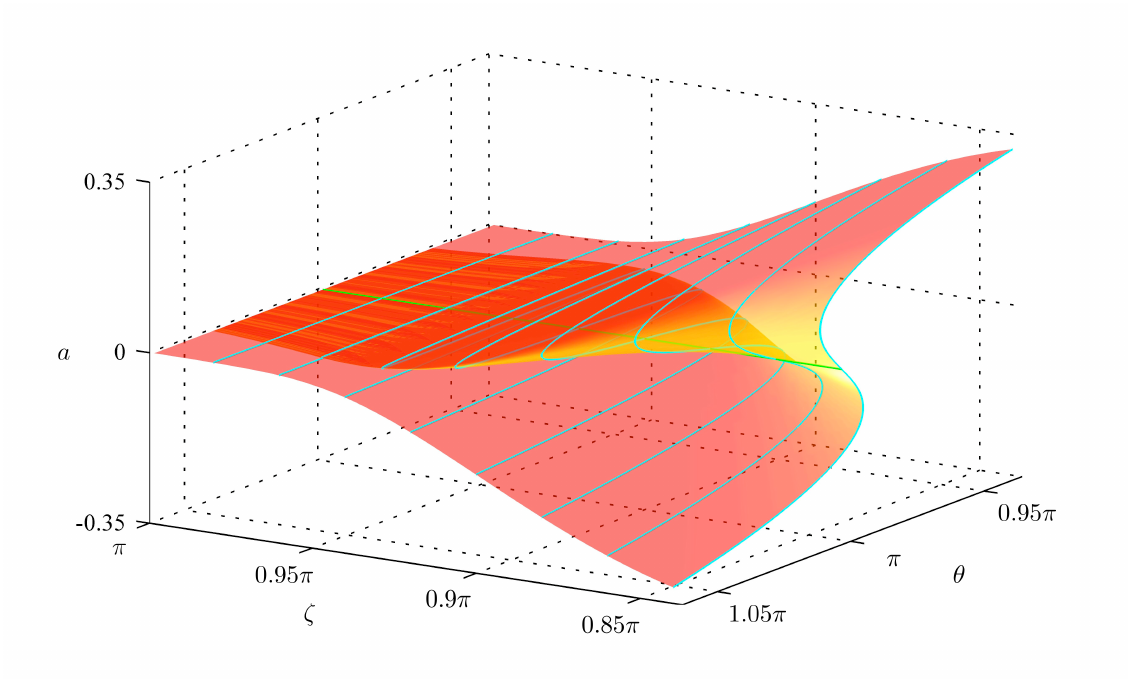} }
    \caption{$\Lambda_{\beta}$ with $\beta = 1$. The black arrows denote unit vectors in the fibers of $\mathcal{N}_{1}$ at various locations along the orbit $E$, which is shown in green. Notice the rotation of the bundle $\mathcal{N}_{1}$ shown in (a). A zoom view of the fold on $\Lambda_1$ is shown in (b) with various cross-sections of $\Lambda_1$ highlighted in light blue.}
  \label{fig:FOLD}
\end{figure}

The above analysis implies the following picture about the rotation of the normal bundle of $E$ under the linearized flow $D \phi_{\beta}^t (0, \zeta, \pi)$. For every $\beta \in [0, 1]$, we obtain a subbundle $\mathcal{N}_{\beta}$ of the normal bundle of $E$ with its fibers oriented according to the values of $\alpha$ along the orbit whose backward-time asymptotic limit is the fixed point $(\alpha, \zeta) = (0, 0)$ in the system (\ref{RT}) with the corresponding $\beta$. Clearly, for any $\beta \in [0, 1]$, the bundle $\mathcal{N}_{\beta}$ is invariant under $D \phi_{\beta}^t (0, \zeta, \pi)$, and its fibers satisfy
\begin{equation*}
\lim_{\zeta \rightarrow 0^+} \mathcal{N}_{\beta}|_{(0, \zeta, \pi)} = \big\{ \lambda \tfrac{\partial}{\partial \theta} : \lambda \in \RR \big\} \,,
\end{equation*}
which coincides with the tangent space of $\Gamma$ at $(a, \zeta, \theta) = (0, 0, \pi)$. Furthermore, for any $\beta \in [0, 1]$, $\mathcal{N}_{\beta}$ is tangent to $\Lambda_{\beta}$ along the orbit $E$. For any $\beta \in [0, \beta_c)$, the bundle $\mathcal{N}_{\beta}$ corresponds to the homoclinic orbit of the fixed point $(\alpha, \zeta) = (0, 0)$ in (\ref{RT}), and thus
\begin{equation*}
\lim_{\zeta \rightarrow 0^+} \mathcal{N}_{\beta}|_{(0, \zeta, \pi)} = \lim_{\zeta \rightarrow \pi^-} \mathcal{N}_{\beta}|_{(0, \zeta, \pi)} \,.
\end{equation*}
In this case, $\Lambda_{\beta}$ is the $C^1$ invariant torus $T_{\beta}$ (see Figure \ref{fig:smooth} for the case with $\beta = 0.65$). However, when $\beta = \beta_c \simeq 0.815$, the bundle $\mathcal{N}_{\beta}$ corresponds to the heteroclinic orbit shown in Figure \ref{fig:RT2} in mixed blue and red and thus
\begin{equation*}
\lim_{\zeta \rightarrow \pi^-} \mathcal{N}_{\beta_c}|_{(0, \zeta, \pi)} = \big\{ \lambda ( - \tfrac{\partial}{\partial a} + \tfrac{\partial}{\partial \theta}) : \lambda \in \RR \big\} \,.
\end{equation*}
Consequently, $\Lambda_{\beta_c}$ loses smoothness at $(a, \zeta, \theta) = (0, 0, \pi)$ and becomes a nonsmooth topological torus (see Figure \ref{fig:nonsmooth}). Finally, for any $\beta \in (\beta_c, 1]$, since the number of rotation in $\alpha$ is $-1$ along the corresponding homoclinic orbit of $(\alpha, \zeta) = (0, 0)$, the bundle $\mathcal{N}_{\beta}$ has a $-\pi$ rotation along the orbit $E$, which causes $\Lambda_{\beta}$ to fold in a neighborhood of $(a, \theta) = (0, \pi)$ in the limit $\zeta \rightarrow \pi^-$ (see Figure \ref{fig:FOLD} for the case with $\beta = 1$).

\begin{thebibliography}{1}

\bibitem{Ch06}
{\sc C.~Chicone}, {\em Ordinary Differential Equations with Applications},
  vol.~34 of Texts in Applied Mathematics, Springer, New York, 2nd~ed., 2006.

\bibitem{ChLi00}
{\sc C.~Chicone and W.~Liu}, {\em On the continuation of an invariant torus in
  a family with rapid oscillations}, SIAM J. Math. Anal., 31 (1999/00),
  pp.~386--415 (electronic).

\bibitem{Fe71}
{\sc N.~Fenichel}, {\em Persistence and smoothness of invariant manifolds for
  flows}, Indiana Univ. Math. J., 21 (1971/1972), pp.~193--226.

\bibitem{HaLl07}
{\sc A.~Haro and R.~de~la Llave}, {\em A parameterization method for the
  computation of invariant tori and their whiskers in quasi-periodic maps:
  explorations and mechanisms for the breakdown of hyperbolicity}, SIAM J.
  Appl. Dyn. Syst., 6 (2007), pp.~142--207 (electronic).

\bibitem{HiPuSh77}
{\sc M.~W. Hirsch, C.~C. Pugh, and M.~Shub}, {\em Invariant Manifolds},
  vol.~583 of Lecture Notes in Mathematics, Springer-Verlag, Berlin, 1977.

\bibitem{Wi94}
{\sc S.~Wiggins}, {\em Normally Hyperbolic Invariant Manifolds in Dynamical
  Systems}, vol.~105 of Applied Mathematical Sciences, Springer-Verlag, New
  York, 1994.

\end{thebibliography}
\end{document}